\documentclass{elsart}
\usepackage{amsmath,amsfonts}
\usepackage[square]{natbib}
\journal{Physics Letters A}

\newcommand\mvector{\boldsymbol}

\newcommand\vd{\mvector{d}}

\newcommand\vp{\mvector{p}}
\newcommand\vq{\mvector{q}}

\newcommand\vv{\mvector{v}}

\newcommand\vx{\mvector{x}}

\newcommand\valpha{\mvector{\alpha}}
\newcommand\vgamma{\mvector{\gamma}}

\newcommand\vlambda{\mvector{\lambda}}
\newcommand\field{\mathbb}
\newcommand\R{\field{R}}

\newcommand\C{\field{C}}
\newcommand\Z{\field{Z}}
\newcommand\N{\field{N}}
\newcommand\Q{\field{Q}}

\newcommand\balpha{\boldsymbol{\alpha}}

\newcommand\bvarphi{\boldsymbol{\varphi}}

\newcommand\rmd{\mathrm{d}}

\newcommand\rmi{\mathrm{i}}

\newcommand\Dt{\frac{\rmd\phantom{t} }{\rmd t}}

\newcommand\gdeg{\operatorname{\mbox{$\vgamma$}-deg}}
 \newcommand\bfi[1]{{\bfseries\itshape{#1\/}}}
\theoremstyle{plain}
\newtheorem{theorem}{Theorem}
\newtheorem{lemma}{Lemma}
\newtheorem{proposition}{Proposition}
\theoremstyle{definition}
\newtheorem{example}{Example}

\begin{document}
\hfill \textsf{accepted to: Physics Letters A}
\begin{frontmatter}
  \title{ Darboux Polynomials and First Integrals of Natural Polynomial
    Hamiltonian Systems }
\author{Andrzej J. Maciejewski}
\address{Institute of Astronomy,
University of Zielona G\'ora \\
Podg\'orna 50, 65--246 Zielona G\'ora, Poland,}
\ead{maciejka@astro.ia.uz.zgora.pl}
\author{Maria Przybylska}
\ead{mprzyb@astri.uni.torun.pl}
\address{INRIA -- Projet {\sc  Caf\'e}\\
2004, Route des Lucioles, B.P. 93\\
06902 Sophia Antipolis Cedex, France\\
and\\
Toru\'n Centre for Astronomy, 
N. Copernicus University,\\
Gagarina 11, 87--100 Toru\'n, Poland}
\begin{abstract}
  We show that for a natural polynomial Hamiltonian system the
  existence of a single Darboux polynomial (a partial polynomial first
  integral) is equivalent to the existence of an additional first
  integral functionally independent with the Hamiltonian function.
  Moreover, we show that, in a case when the degree of potential is odd,
   the system does not admit any proper Darboux polynomial, i.e.,
  the only Darboux polynomials are first integrals.
\end{abstract}
\end{frontmatter}
\section{Introductions}
\label{sec:intro}
Let us consider a system of ordinary differential equations
\begin{equation}
\label{eq:ds}
\Dt \vx = \vv(\vx),  \qquad  \vx=(x_1,\ldots, x_n)\in U\subset\C^n, 
\end{equation}
where the right hand sides $\vv(\vx)=(v_i(\vx),\ldots,v_n(\vx))$ are
smooth on a domain $U$ in $\C^n$.  Let $\bvarphi: \R\supset I \mapsto
U$ be a solution of~\eqref{eq:ds} defined on a certain open non-empty
interval $I$ of the real axis.  A continuous function $F:U\rightarrow \R$
is called a \bfi{first integral} of system~\eqref{eq:ds} if it is
constant along its solutions, i.e., if function $ F\circ\bvarphi$ is
constant on its domain of definition for an arbitrary solution $\bvarphi$
of~\eqref{eq:ds}.  When  $F$ is differentiable, then it is
a first integral of system~\eqref{eq:ds} if
\begin{equation}
\label{eq:Lv}
L_{\vv}(F)(\vx):= \sum_{i=1}^n v_i(\vx)\partial_iF(\vx)=0,  \qquad \partial_i \equiv \frac{\partial\phantom{x}}{\partial x_i}, 
\end{equation}
where $L_{\vv}$ is the Lie derivative along vector field $\vv$. 
If $F$ is a first integral of~\eqref{eq:ds}, then its constant value levels
\begin{equation}
M_c(F)=\{ \vx\in U\, |\quad F(\vx)=c \,\} ,  
\end{equation}
are invariant with respect to the flow generated by~\eqref{eq:ds}. It
appears that together with the first integrals it is worth looking for
partial first integrals. We say that $F$ is a \bfi{partial first
  integral} if its zero level $M_0(F)$ is invariant with respect to
the flow generated by~\eqref{eq:ds}. A differentiable function $F$ is
a partial first integral of system~\eqref{eq:ds} if $L_{\vv}(F)(\vx)=
0$ for all $\vx$ such that $ F(\vx)=0$.

There are a lot of practical and theoretical reasons why we would like
to know if a given system of differential equations possesses a first
integral. However, in its most complete and general case this problem is
intractable. We can simplify it a little bit, e.g., by imposing
additional conditions on the smoothness of the right hand sides
of~\eqref{eq:ds} and smoothness of the first integrals. Nevertheless, even
if we assume that the right hand sides of the investigated system are
polynomial and we look for a polynomial first integral the problem is
very hard to approach.

In this paper we investigate systems~\eqref{eq:ds} of a special form
with the polynomial right hand sides. Namely, we consider polynomial
Hamiltonian systems defined on $\C^{2m}$. Let $\vq=(q_1,\ldots,q_m)$
and $\vp=(p_1,\ldots,p_m)$ denote the standard canonical coordinates
and $H=H(\vq,\vp)$ is a polynomial Hamiltonian function.  We say that
a first integral $F=F(\vq,\vp)$ of the canonical Hamiltonian equations
is an \bfi{additional first integral} if $F$ and $H$ are functionally
independent.

It was observed \cite{Kozlov:96::,Perelomov:90::} that for most known
polynomial Hamiltonian systems an additional first integral is a
polynomial one and its degree with respect to $\vp$ is small. However,
there are also examples of polynomial Hamiltonian systems for which an
additional first integral is more complicated.
\begin{example}
For the Hamiltonian system given by 
\[
H = \frac{1}{2}(p_1^2+p_2^2) + 2q_2p_1p_2 - q_1
\]
an additional first integral has the following form
\[
 F = p_2\exp(p_1^2),
\]
see \cite{Perelomov:90::}.
\end{example}
\begin{example}
For the Hamiltonian system given by
\[
H = \frac{1}{2}(p_1^2+p_2^2) -f(p_1,p_2)(p_1q_1 -\alpha p_2q_2), \qquad 
\alpha\in\C, 
\]
where $f(p_1,p_2)$ is an arbitrary polynomial, 
an additional first integral has the following form
\[
 F = p_1^{\alpha}p_2.
\]
If $\alpha\not\in\Q$ this first integral is not meromorphic, when
$\alpha\in\Q$, then it is algebraic, when $\alpha\in-\N$, then it is
rational, and finally when $\alpha\in\N$, then it is a polynomial of
degree $\alpha+1$ with respect to the  momenta. 
\end{example}
These examples show that it is difficult to expect that the additional
first integrals of a polynomial Hamiltonian system should be distinguished
from the first integrals of an arbitrary polynomial system.

In mechanics there exists a well defined class of natural Hamiltonian
systems for which Hamilton's function is a sum of the kinetic and
potential energy. In this paper we consider systems for which  the 
Hamiltonian has the following form
\begin{equation}
  \label{eq:nham}
  H = \frac{1}{2}\sum_{i=1}^m\mu_i p_i^2 + V(\vq), 
\end{equation}
where $V(\vq)$ is a polynomial, and $\mu_i\in\C$ for $i=1,\ldots, m$.
To avoid trivialities we assume that $m>1$ and $\deg V(\vq)>2$. For
such systems canonical equations have the form
\begin{equation}
\label{eq:heqs}
\begin{split}
  \Dt q_i &=  \mu_i p_i, \\
  \Dt p_i &= -\frac{\partial V}{\partial q_i}, \qquad i=1, \ldots, m.
\end{split}
\end{equation}

We did our best but we have not find even a single example of a system
with the Hamiltonian of the form~\eqref{eq:nham} admitting an additional
first integral which is not a polynomial one. Thus, it seems 
reasonable to formulate the following questions:
\begin{enumerate}
\item Does there exist a natural polynomial Hamiltonian system
  admitting a proper rational additional first integral?
\item Does the existence of a proper rational additional first
  integral of a natural polynomial Hamiltonian system imply the existence 
  of an additional  polynomial first integral?
\item Does there exist a natural polynomial Hamiltonian system admitting an
  analytic additional first integral which is not a polynomial one?
\item Does the existence of an analytic additional first integral of a
  natural polynomial Hamiltonian system imply the existence of  a
  polynomial one?
\item Does there exist a natural polynomial Hamiltonian system admitting a
  meromorphic  additional first integral which is not a rational one?
\item Does the existence of a meromorphic additional first integral of a
  natural polynomial Hamiltonian system imply the existence of a
  rational one?
\end{enumerate}
By a proper rational first integral we mean a rational first integral
which is not a ratio of two polynomial first integrals. In the above
questions we do not mention  algebraic first integrals, because, as it
is well known, the existence of an algebraic additional first integral 
 implies the existence of a rational one, see e.g. \cite{Churchill:90::c}. 
  
In this paper we focused our attention on the first two questions.
Roughly speaking, the main result of this paper is an affirmative
answer to the second question. 

In fact our result is stronger. Let us denote by $L_H$ the Lie
derivative along the Hamiltonian vector field defined
by~\eqref{eq:heqs}.  A polynomial $G=G(\vq,\vp)\neq 0$ is called a
\bfi{Darboux polynomial} of system~\eqref{eq:heqs} if
\begin{equation}
  \label{eq:dlh}
   L_H G = \Lambda G
\end{equation}
for a certain polynomial $\Lambda=\Lambda(\vq,\vp)$, which is called
\bfi{cofactor} of $G$. When $\Lambda=0$, then the Darboux polynomial
is a first integral. We say that $G$ is a \bfi{proper Darboux
  polynomial} if $\Lambda\neq0$.  Let us assume that polynomials
$G=G(\vq,\vp)$ and $Q=Q(\vq,\vp)$ are coprime, i.e., $G$ and $Q$ do
not have common divisors different from constants. Then it is easy to show
that if $F=G/Q$ is a proper rational first integral of
system~\eqref{eq:heqs}, then $G$ and $Q$ are proper Darboux
polynomials with the same cofactor $\Lambda$.

The main results of this paper are formulated in the following theorems.
\begin{theorem}
\label{thm:t1}
Assume that in Hamiltonian~\eqref{eq:nham} potential $V(\vq)$ is of 
odd degree, then every Darboux polynomial of Hamiltonian
system~\eqref{eq:heqs} is a first integral.
\end{theorem}
\begin{theorem}
\label{thm:t2}
Assume that in Hamiltonian~\eqref{eq:nham} at least two $\mu_i$ are
not zero and potential $V(\vq)$ and is of even degree. If Hamiltonian
system~\eqref{eq:heqs} possesses a proper Darboux polynomial, then it
admits an additional polynomial first integral.
\end{theorem}
The statements of the above theorems are rather amazing. In fact, it
is obvious that the conditions for the existence of a Darboux
polynomial are much weaker than the conditions for the existence of a
polynomial first integral. The first theorem says that it is not the
case for natural polynomial Hamiltonian systems with a potential of
odd degree.  Moreover, let us notice that in the second theorem we do
not exclude the existence of a proper Darboux polynomial (there are
many examples where it exists). However, for the existence of a proper
rational first integral we need \bfi{two} proper Darboux polynomials
with the \bfi{same} cofactor. The second theorem says that one Darboux
polynomial is enough to construct an additional polynomial first
integral.

Proofs of the above theorems are purely algebraic and this is  why we
collect all necessary facts from (elementary) differential algebra in
the next section. In Section~3 we give proofs of our theorems. In the
last section we show some generalizations of our results and we give
several remarks.

\section{Basic notions of differential algebra}
\label{sec:basic}
When the algebraic methods are used to study a polynomial or rational
differential equation, it is convenient to use tools of differential
algebra. It gives a natural language for the problem, simplifies
exposition and proofs.  In this section we fix notation and we collect
basic facts from differential algebra. We follow the excellent
lectures of A.~Nowicki\footnote{The whole book can be found at
  \texttt{http://www.mat.uni.torun.pl/\textasciitilde anow/polder.html
  }}\cite{Nowicki:94::} simplifying only some assumptions. All
propositions formulated in this section are very simple. We left their
proofs to the readers, and we refer them to~\cite{Nowicki:94::} and
references therein.
\subsection{Notation}
By $\N_0$ we denote the set of non-negative integers.  If $k$ is  a
field (we always assume that $k=\R$ or $k=\C$) then $ k[x_1,\ldots,
x_n]$ and $ k(x_1,\ldots, x_n )$ denote the ring of polynomials and the 
field of rational functions, respectively. We use also an abbreviated
notation, i.e., $k[\vx]= k[x_1,\ldots, x_n]$ and $k(\vx)=
k(x_1,\ldots, x_n )$.

If $\balpha=(\alpha_1,\ldots,\alpha_n)\in\N_0^n$, then we denote
\[
\vx^{\balpha}=x_1^{\alpha_1}\cdots x_n^{\alpha_n}, \qquad
|\balpha|=\sum_{i=1}^n \alpha_i.
\] 
A non-zero element $\vgamma\in\N^n $ we called a \bfi{direction}, and
we denote
\[
\langle \vgamma, \balpha\rangle= \sum_{i=1}^n \gamma_i\alpha_i,
\qquad\text{for}\quad\balpha\in\N_0^n.
\]
A non-zero polynomial $F\in k[\vx]$ is called a \bfi{$\vgamma$-form of
  $\vgamma$-degree $r$} if it has the form
\[
F= \sum_{\langle \vgamma, \valpha\rangle=r}F_{\valpha} \vx^{\valpha},
\qquad F_{\valpha}\in k.
\]
In such a case we also say that $F$ is a \bfi{$\vgamma$-homogeneous
  polynomial of $\vgamma$-degree $r$}.  The zero polynomial is a
$\vgamma$-form of an arbitrary degree.
\begin{proposition}
  If $F$ is a non-zero polynomial in $k[\vx]$, then the following
  conditions are equivalent:
\begin{enumerate}
\item $F$ is a $\vgamma$-form of $\vgamma$-degree $r$.
\item $F(t^{\gamma_1}x_1, \ldots, t^{\gamma_n}x_n)= t^r
  F(x_1,\ldots,x_n)$.
\item $ \gamma_1x_1\partial_1F+\cdots +\gamma_nx_n\partial_nF=rF$.
\end{enumerate}
\end{proposition}
The introduced $\vgamma$-homogeneity converts $k[\vx]$ into a
\bfi{$\vgamma$-graded ring}.  Namely, we have
\[
k[\vx]=\bigoplus_{r\in\N_0} k_r^{\vgamma}[\vx],
\]
where $k_r^{\vgamma}[\vx]$ denotes the group all $\vgamma$-forms of
$\vgamma$-degree $r$, and
\[ 
k_r^{\vgamma}[\vx] k_s^{\vgamma}[\vx]\subseteq k_{r+s}^{\vgamma}[\vx],
\]
for all $r,s\in\N_0$. If there is no ambiguity, we write $k_r[\vx]$
instead of $k_r^{\vgamma}[\vx]$.  Thus, every polynomial $F\in k[\vx]$
has the following \bfi{$\vgamma$-decomposition}
\[
F = \sum_{s\in S} F_s,
\]
where $F_s$ is a $\vgamma$-form of $\vgamma$-degree $s$, and $S$ is a
finite set of integers.  By $\gdeg(F)$ we denote the \bfi{
  $\vgamma$-degree} of $F$, that is, for $F\neq0$, the maximal $s$
such that $F_s\neq 0$ in the above $\vgamma$-decomposition of $F$; for
$F=0$ we put $\gdeg(F)=-\infty$. If $s= \gdeg(F)$ then we denote
$F_s=F^+$.

If $F,G\in k[\vx]$ and $FG$ is a non-zero $\vgamma$-form, then $F$
  and $G$ are $\vgamma$-forms. 

  Let us fix two $\vgamma$-forms $F$ and $G$ of $\vgamma$-degree $p$
  and $q$, respectively.  In the ring $k[X,Y]$ we define
  $\vlambda$-gradation putting $\vlambda=(p,q)$. For a polynomial
  $W\in k[X,Y]$ we define a polynomial $\widetilde W\in k[\vx]$ as
  $\widetilde W=W[F,G]$.
\begin{proposition}
\label{prop:gl}
Let $P\in k[X,Y]$, and
\[
  P = \sum_{s\in S} P_s, 
\]
where $S$ is a finite set of integers, be its   $\vlambda$-decomposition. 
If $\widetilde P =0 $, then  $\widetilde P_s= 0$ for $s\in S$.
\end{proposition}
\begin{pf}
  For each $i,j\in\N_0$ polynomials $F^iG^j$ are $\vgamma$-forms in $k[\vx ]$
  of $\vgamma$-degree $ip+jq$.  Thus, a $\vgamma$-homogeneous component
  of $\widetilde P$ of $\vgamma$-degree $s$ is just $ \widetilde
  P_s$.\qed
\end{pf}   

We say that rational functions $F_i \in k(\vx)$, $i=1, \ldots, m$ are
\bfi{functionally independent} if their differentials $\rmd F_i$ are
linearly independent over $k$, i.e., if matrix
\[
\begin{bmatrix}
\partial_1 F_1 & \ldots & \partial_n F_1 \\
\ldots & \ldots & \ldots \\
\partial_1 F_m & \ldots & \partial_n F_m \\
\end{bmatrix} 
\]
has the maximal rank. 

We say that rational functions $F_i \in k[\vx]$, $i=1, \ldots, m$ are
\bfi{algebraically dependent over $k$} if there exists a non-zero polynomial 
$W\in k[X_1,\ldots, X_m]$ such that $W(F_1(\vx), \ldots, F_m(\vx)) =0$. 

Later, we use the following lemma, see e.g. Proposition~1.16 in \cite{Churchill:96::b}. 
\begin{lemma}
In $\C(\vx)$ functional and algebraic independences are equivalent.  
\end{lemma}
\subsection{Derivations.}

Let us fix $n$ polynomials $v_i\in k[\vx]$, $i=1,\ldots,n$, and denote
$\vv=(v_1,\ldots,v_n)$. Then, let $\vd_{\vv}: k[\vx]\mapsto k[\vx]$ be
a map defined in the following way
\begin{equation}
\vd_{\vv}(F) = \sum_{i=1}^n v_i \partial_i F, 
\qquad\text{for}\quad F\in  k[\vx].
\end{equation}
Obviously it is $k$-linear and satisfies the Leibnitz rule, i.e., we
have
\[
\vd_{\vv}(FG) = \vd_{\vv}(F)G + F\vd_{\vv}(G),
\] 
for all $ F,G\in k[\vx]$. Such maps are called derivations. For a more
general definition see \cite{Nowicki:94::}.  In the ring $k[\vx]$
every derivation is of this form.
For an element $R\in k(\vx)$ given by $R=F/G$, where $F$ and $G$ are
coprime polynomials, we put
\[
\vd_{\vv}(F) = \frac{ \vd_{\vv}(F)G - F\vd_{\vv}(G)}{R^2}.
\]
In this way we extend uniquely $\vd_{\vv}$ into $k(\vx)$.  If it is
clear from the context,  we write $\vd$ instead of $\vd_{\vv}$.  The
kernel of $\vd$ in $k[\vx]$, which we denote
\[
k[\vx]^{\vd} = \{ F\in k[\vx] \,|\, \vd(F)=0 \},
\]
is a subring of $k[\vx]$.  It is called the \bfi{ring of constants of
  derivation} $\vd$.  Similarly, the kernel of $\vd$ in $k(\vx)$,
which we denote
\[
k(\vx)^{\vd} = \{ R\in k(\vx) \,|\, \vd(R)=0 \},
\]
is a subfield of $k(\vx)$, and it is called the \bfi{field of constants of
  derivation} $\vd$.

As we see, having a system of polynomial equations~\eqref{eq:ds} with
right hand sides $\vv$, we associate with this system a derivation
$\vd=\vd_{\vv}$ which is, in this context,  nothing else but the Lie
derivative along $\vv$ considered as a differential operator.  Thus, $
k[\vx]^{\vd}$ and $k(\vx)^{\vd}$ are  sets of the polynomial and rational
first integrals, respectively.

Now, the notion of Darboux polynomial is defined in the following way.
We say that non-zero polynomial $F\in k[\vx]$ is a \bfi{Darboux
  polynomial} of derivation $\vd$ if $\vd(F)=\Lambda F$ for a certain
$\Lambda\in k[\vx]$. Polynomial $\Lambda$ is called a \bfi{cofactor}
of Darboux polynomial $F$.

The following proposition explains the importance of the Darboux polynomials.
\begin{proposition}
\label{prop:rat}
  Let $F$ and $G$ be non-zero coprime polynomials in $k[\vx]$. Then
\[
\vd(F/G) = 0,
\]
if and only if
\[
\vd(F) = P F \quad\text{and}\quad \vd(G) = P G,
\]
for some $P\in k[\vx]$.
\end{proposition}
Two important properties of Darboux polynomials are described in the
following propositions.
\begin{proposition}
\label{prop:prod}
 If $\vd (F_i) = \Lambda_i F_i$ for $i=1,2$, then $\vd (F_1F_2) =
  ( \Lambda_1 + \Lambda_2)(F_1F_2)$.
\end{proposition}
\begin{proposition}
\label{prop:fact}
If $F\in k[\vx]$ is a Darboux polynomial of derivation $\vd$, then all
factors of $F$ are also  Darboux polynomials of derivation $\vd$.
\end{proposition}
Now, we fix a direction $\vgamma$. We say that derivation $\vd$ is
\bfi{$\vgamma$-homogeneous of $\vgamma$-degree $s$} if
\[
\vd( k_r^{\vgamma}[\vx] ) \subseteq k_{r+s}^{\vgamma}[\vx],
\]
for an arbitrary $r\in\Z$.  For checking if a given derivation is
$\vgamma$-homogeneous it is convenient to apply the following proposition.
\begin{proposition}
 \label{prop:dhom}
 Derivation $\vd$ is $\vgamma$-homogeneous of $\vgamma$-degree $s$ if
 and only if $\vd(x_i)$ is a $\vgamma$-form of $\vgamma$-degree
 $s+\gamma_i$, for $i=1,\ldots,n.$
\end{proposition}
Homogeneity of a derivation is important when we look for its Darboux
polynomials. 
\begin{proposition}
 \label{prop:hom}
 If $F\in k[\vx]$ is a Darboux polynomial of a $\vgamma$-homogeneous
 derivation $\vd$ of $\vgamma$-degree $s$, then a cofactor $\Lambda\in
 k[\vx]$ is $\vgamma$-form of $\vgamma$-degree $s$, and all
 $\vgamma$-components of $F$ are also Darboux polynomials with the
 common cofactor $\Lambda$.
\end{proposition}
For a polynomial derivation $\vd_{\vv}$ and a given $\vgamma$-gradation
we can always represent $\vd_{\vv}$ as a finite sum of $\vgamma$
homogeneous derivations. Namely, for an integer $l$ we denote 
\[
\vv^{(\vgamma,l)}=(v^{(\vgamma,l)}_1, \ldots,v^{(\vgamma,l)}_n),
\] 
where $v^{(\vgamma,l)}_i$ is the $\vgamma$-homogeneous component of
$v_i$ of $\vgamma$-degree $l+\gamma_i$. Then by
Proposition~\ref{prop:dhom} $\vd_{\vv^{(\vgamma,l)}}$ is a
$\vgamma$-homogeneous derivation of $\vgamma$-degree $l$.  By
$\vd_{\vv}^+$ we denote a $\vgamma$-homogeneous derivation
$\vd_{\vv^{(\vgamma,l)}}$ corresponding to the maximal value of $l$.
By $\vgamma$-degree of $\vd_{\vv}$ we understand $\vgamma$-degree of
$\vd_{\vv}^+$
\begin{proposition}
\label{prop:ho}
Assume that $F\in k[\vx]$ is a Darboux polynomial of a derivation
$\vd_{\vv}$ and $\Lambda\in k[\vx]$ is a corresponding cofactor. Then for
an arbitrary $\vgamma$-gradation we have $\vd_{\vv}^+(F^+)=P^+F^+$.
\end{proposition}
\begin{proposition}
\label{pro:degl}
  If $\vd_{\vv}(F)=\Lambda F$, then $\gdeg\Lambda\leq r$ where $r$ is
  $\vgamma$-degree of $\vd_{\vv}$. 
\end{proposition}
\section{Proofs}
To translate the language of Hamiltonian systems into the framework of
differential algebra, we introduce the following conventions. By
$k[\vq,\vp]$ we denote the ring of polynomials of $n=2m$ variables
\[
\vx=(\vq,\vp)=(q_1,\dots q_m,p_1,\ldots, p_m).
\]
For a polynomial Hamiltonian $H \in k[\vq,\vp]$ we define a derivation
$\vd_H$ in the following way
\begin{equation*}
\label{eq:d_H}
\vd_{H}(q_i) = \frac{\partial H}{\partial p_i}, \qquad 
\vd_{H}(p_i) = -\frac{\partial H}{\partial q_i}, \qquad i=1,\ldots, m. 
\end{equation*}  
Thus,  for Hamilton's function~\eqref{eq:nham} the explicit form of
 $\vd_H$ is following
\begin{equation}
\label{eq:d_Hnat}
\vd_{H}(q_i) =  \mu_i p_i, \qquad
\vd_{H}(p_i) = -\frac{\partial V}{\partial q_i}, \qquad i=1,\ldots, m. 
\end{equation}  
Natural Hamiltonian systems have a very important property of the time
reversibility. We can define this property in the following way.  Let
$\tau: k[\vq,\vp] \mapsto k[\vq,\vp] $ denote an automorphism of
$k[\vq,\vp] $ given by:
\[
\tau(q_i) = q_i, \qquad \tau(p_i)=-p_i, \qquad i=1, \ldots, m.
\]
The following proposition is easy to check.
\begin{proposition}
\label{prop:tr}
Derivation~\eqref{eq:d_Hnat} has  the following property
\[
\tau^{-1} \circ \vd_{H} \circ \tau = - \vd_{H}.
\]
\end{proposition}
For a given  Hamilton's function~\eqref{eq:nham}   we fix 
a $\vgamma$-gradation defined by the following direction
\[
\vgamma=(\gamma_1, \ldots, \gamma_m, \gamma_{m+1}, \ldots, \gamma_{2m})=
(2,\ldots,2,r, \ldots,r), 
\]
where $r=\gdeg V$. 
We prove the following simple but important lemma.
\begin{lemma}
  Assume that potential $V\in k[\vq]$ is homogeneous of degree $r>2$.
  Then derivation $\vd_H$ given by \eqref{eq:d_Hnat} is
  $\vgamma$-homogeneous of degree $s=r-2$.
\end{lemma}
\begin{pf}
  Variables $q_i$ and $p_i$ are $\vgamma$-forms of degree $2$ and $r$,
  respectively. Thus,
\[
\gdeg (\vd_H(q_i)) = \gdeg (p_i) = r=s+\gamma_i,
\]
as $\gamma_i=2$ for $i=1,\ldots, m$. Moreover, we have
 \[
 \gdeg( (\vd_H(p_i) )=\gdeg\left( -\frac{\partial V}{\partial q_i}
 \right) = 2(r-1)= (r-2) + r= s + \gamma_{i+m},
 \]
 as $\gamma_{i+m}=r$ for $i=1,\ldots,m$. Now, the proof follows from
 Proposition~\ref{prop:dhom}. \qed
\end{pf}
It is easy to show that
\[
H^+= \frac{1}{2}\sum_{i=1}^m \mu_i p_i^2 + V^+(\vq),
\]
where $V^+(\vq)\neq0$ is the homogeneous component of $V(\vq)\neq0$ of
the highest degree $r$. An important consequence of the previous lemma is
following.
\begin{lemma}
\label{lem:P}
  Let $\deg V \geq 2 $ and
  $F\in k[\vq,\vp]$ be a Darboux polynomial of $\vd_H$ with cofactor  $\Lambda$. Then $\Lambda\in
  k[\vq]$.
\end{lemma}
\begin{pf}
  We know that $\vgamma$-degree of $\vd_H$ is $s=r-2$. Thus, by
  Proposition~\ref{pro:degl} $\gdeg \Lambda\leq s=r-2$.  But $\gdeg
  p_i=r>s$ for $i=1,\ldots, m$. Hence, $\Lambda$ does not depend on
  $\vp$.  \qed
\end{pf}
Now, the proof of Theorem~\ref{thm:t1} is elementary.
\begin{pf*}{PROOF of Theorem~\ref{thm:t1}.}
  Let us assume the opposite, i.e., let $F$ be a Darboux polynomial of
  $\vd_H$, with non-zero cofactor $\Lambda$.  Then, by
  Proposition~\ref{prop:ho}, we have $\vd_{H^+}F^+= \Lambda^+F^+$.
  Hence, by Proposition~\ref{prop:hom}, $\Lambda^+$ is a non-zero
  $\vgamma$-form and $\gdeg \Lambda^+ = r-2$, i.e., for an odd $r$ it
  is of odd degree. But, by~Lemma~\ref{lem:P}, $\Lambda^+$ depends
  only on $\vq$, so its $\vgamma$-degree is necessarily even. We have
  the contradiction. Thus, we must have $\Lambda=0$.\qed
\end{pf*}
When the degree of the potential is even, then derivation $ \vd_{H}$ can
possess a proper Darboux polynomial.  However, in such a case  we have 
 also a polynomial first integral.
\begin{lemma}
\label{lem:tau}
  Assume that the degree of $V$ is even and that $ \vd_{H}$ has a
  proper Darboux polynomial $F$ with cofactor $\Lambda$. Then $
  \tau(F)$ is also a proper Darboux polynomial with cofactor
  $-\Lambda$, and $G=\tau(F)F$ is a polynomial first integral of $
  \vd_{H}$.
\end{lemma}
\begin{pf}
If $ \vd_{H}(F) = \Lambda F$, then  we have
\[
\tau(-\vd_{H}(F)) = \tau(-\Lambda F) = - \tau(\Lambda) \tau(F).
\]
But, from Proposition~\ref{prop:tr}, it follows that
\[
\tau(-\vd_{H}(F)) = \tau(\tau^{-1} \circ \vd_{H} \circ \tau (F))=
\vd_{H}(\tau(F)).
\]
Hence $ \vd_{H}(\tau(F)) = - \tau(\Lambda) \tau(F)$ . However, by
Lemma~\ref{lem:P}, cofactor $\Lambda$ does not depend on $\vp$, so
$\tau(\Lambda) =\Lambda$.  As a result we have
\[
\vd_{H}(\tau(F)) = - \Lambda\tau(F).
\]
Now, the fact that $G=\tau(F)F$ is a first integral of $ \vd_{H}$ follows
from Proposition~\ref{prop:prod}.  \qed
\end{pf}
At this point, the question is if the first integral constructed in the above
Lemma is functionally independent of $H$.  To show this we need to
use the fact that the Hamiltonian~\eqref{eq:nham} is irreducible.
\begin{lemma}
\label{lem:irr}
If $m\geq 2$, $V\neq 0$, and among numbers $\mu_i$, $i=1, \ldots, m$,
at least two are non-zero, then  polynomial $H\in k[\vq,\vp]$ given by
~\eqref{eq:nham} is not reducible.
\end{lemma}
\begin{pf}
  Without loss of generality we can assume that $\mu_1\neq 0$ and
  $\mu_2\neq 0$.  Assume that $H$ is reducible. As $H$ is a polynomial
  of second degree with respect to $ \vp$, its factors are polynomials
  of degree one with respect to $\vp$.  So, we can write
\begin{equation}
 2H = \sum_{i=1}^m\mu_i p_i^2 + 2V(\vq)=
 \left(  \sum_{i=1}^m\alpha_i p_i + W_1\right)
 \left(  \sum_{i=1}^m\beta_i p_i +W_2\right),
\end{equation}
where $\alpha_i, \beta_i\in\C$,  for $i=1,\ldots,m$,  $W_1, W_2\in k[\vq]$, and 
$W_1W_2=2V$.  Comparing the left and right hand sides in the above equality we find,
among others, that
\begin{gather}
\label{eq:war1}
\mu_1 = \alpha_1\beta_1 , \qquad   \mu_2 = \alpha_2\beta_2 , \qquad
\alpha_1\beta_2 +\alpha_2\beta_1 = 0 ,\\
\label{eq:war2}
\alpha_1 W_2  + \beta_1 W_1 = 0, \qquad \alpha_2 W_2  + \beta_2W_1 = 0,
\end{gather}
From equations~\eqref{eq:war1} it follows that
\begin{equation}
\label{eq:aimui}
\alpha_1^2 \mu_2 + \alpha_2^2 \mu_1=0.
\end{equation}
Then, multiplying the first equation~\eqref{eq:war2} by
$\alpha_1\mu_2$, and the second one by $\alpha_2\mu_1$ and adding them
together we obtain
\[
(\alpha_1^2 \mu_2 + \alpha_2^2 \mu_1) W_2 + 2\mu_1\mu_2 W_1 = 0.
\]
Hence, taking into account  \eqref{eq:aimui} we have
\[
2\mu_1\mu_2 W_1 = 0.
\]
But, $\mu_1\neq0$, $\mu_2\neq 0$ , and $W_1\neq0$.  The contradiction
finishes the proof. \qed
\end{pf}
\begin{lemma}
\label{lem:rr}
Let us assume that Hamiltonian $H\in k[\vq,\vp]$ satisfies the following
conditions:
\begin{enumerate}
\item $V\neq 0$, is a homogeneous potential of even degree;
\item $m\geq 2$  and among numbers $\mu_i$, $i=1, \ldots, m$
at least two are non-zero; 
\item there exists a proper Darboux  polynomial $F$ of $\vd_H$. 
\end{enumerate} 
Then the first integrals $H$ and $ G=\tau(F)F$ of $\vd_H$ are
algebraically independent.
\end{lemma}
\begin{pf}
  Let us assume that $H$ and $G$ are algebraically dependent and let
  $W\in k[X,Y]$ be the minimal polynomial, i.e. $\widetilde W=
  W[H,G]=0$. From Proposition~\ref{prop:gl} it follows that we can
  assume that $W$ is $\vlambda$-homogeneous. As $W$ is minimal, it is
  irreducible and this is why it can be written in the following form
  $W= X^i + YW_1(X,Y)$ where $i>0$ is an integer and $W_1\neq 0$ is an
  polynomial in $k[X,Y]$. Hence, we have 
\[
 GW_1(H,G) = -H^i.
\]
Thus irreducible factors of $G$ divide $H$. But it is impossible as,
under the assumptions we made, by Lemma~\ref{lem:irr}, $H$ is irreducible.
A contradiction finishes the proof. \qed
\end{pf}
Now, we are ready to prove Theorem~\ref{thm:t2}. 

\begin{pf*}{PROOF of Theorem~\ref{thm:t2}.}
  Let $F$ be a proper Darboux polynomial of $\vd_H$ with cofactor
  $\Lambda$, so $G=\tau(F)F$ is a first integral of $\vd_H$ .  Then,
  by Proposition~\ref{prop:ho}, we have $\vd_{H^+} F^+ =
  \Lambda^+F^+$.  Hence, by Lemma~\ref{lem:tau}, $G^+=\tau(F^+)F^+$ is
  a polynomial first integral of $\vd_{H^+}$. Moreover, by
  Lemma~\ref{lem:irr}, $G^+$ and $H^+$ are algebraically, and thus
  functionally  independent.  But this implies that $G$ and $H$ are
  functionally independent.
\end{pf*} 
\section{Remarks and comments}
There is an important consequence of Theorem~\ref{thm:t1}.  Namely, if
a system with Hamilton's function given by~\eqref{eq:nham} and the
potential of odd degree possesses an additional first integral then it
possesses an additional irreducible first integral. In fact, for such
system, by Theorem~\ref{thm:t1}, all factors of the a polynomial first
integral are first integrals.

In a case of the potential of even degree, an additional first integral can
be reducible or irreducible. As an example, let us consider a Hamiltonian
system with two degrees of freedom with a homogeneous potential of
degree 4
\begin{equation}
H = \frac{1}{2}(p_1^2 + p_2^2) + V(q_1,q_2).
\end{equation}

The list of integrable systems of this form with a first integral of
degree at most 4 in momenta is given in \cite{Hietarinta:83::}. We
have the following possibilities.
\begin{enumerate}
\item When $V= q_1^4$ , then the additional first integral $F=p_2$ is
  irreducible.
\item When $V = (q_1^2+q_2^2)^2$, then the additional first integral
  $F=q_1p_2-q_2p_1$ is irreducible.
\item If $V = C q_1^2 + q_2^4$, then the additional reducible first integral
  $F= p_2^2 + 2q_2^2 = (\rmi p_2 +\sqrt{2}q_2^2)(-\rmi p_2
  +\sqrt{2}q_2^2)$. Obviously in this case $G_1 = \rmi p_2
  +\sqrt{2}q_2^2$ and $G_2 = \tau(G_1)$ are proper
  Darboux polynomials.
\item For $ V = 4q_1^4/3 + q_1^2q_2^2 + q_2^4/12$  the
  additional first integral $F= p_2(p_1q_2-p_2q_1)
  +q_2^2(2q_1^3+q_1q_2^2)/3$ is irreducible.
\item When $ V = 4q_1^4/3 + q_1^2q_2^2 + q_2^4/6$ the additional first
  integral $F$ of degree 4 with respect to momenta is reducible and
  can be written in the form $F = \tau(G) G$ where
 \[
G = 3\sqrt{6}p_2^2 + 12\rmi p_2q_1q_2 + q_2^2\left[ -6\rmi p_1+\sqrt{6}(2q_1^2+q_2^2)\right]. 
\]
\end{enumerate}
\begin{ack}
As usual, we thank Zbroja not only for her linguistic help. For the
second author this research has been supported by a Marie Curie
Fellowship of the European Community programme Human Potential under
contract number HPMF-CT-2002-02031.
\end{ack}
%
\newcommand{\noopsort}[1]{}\def\cprime{$'$}
  \def\cydot{\leavevmode\raise.4ex\hbox{.}}

\end{document}